\documentclass[a4paper,notitlepage,11pt]{article}%
\usepackage{amssymb}
\usepackage{graphicx}
\usepackage{amsmath}
\usepackage{amsthm}
\usepackage{amsfonts}%
\providecommand{\U}[1]{\protect\rule{.1in}{.1in}}
\theoremstyle{plain}
\newtheorem{theorem}{Theorem}[section]

\newtheorem{corollary}[theorem]{Corollary}

\newtheorem{lemma}[theorem]{Lemma}

\theoremstyle{definition}

\newtheorem{remark}[theorem]{Remark}
\numberwithin{equation}{section}
\numberwithin{theorem}{section}
\ifx\pdfoutput\relax\let\pdfoutput=\undefined\fi
\newcount\msipdfoutput
\ifx\pdfoutput\undefined\else
\ifcase\pdfoutput\else
\ifx\paperwidth\undefined\else
\ifdim\paperheight=0pt\relax\else\pdfpageheight\paperheight\fi
\ifdim\paperwidth=0pt\relax\else\pdfpagewidth\paperwidth\fi
\fi\fi\fi
\begin{document}

\title{Strictly positive solutions for one-dimensional nonlinear elliptic problems
\thanks{2000 \textit{Mathematics Subject Clasification}. 34B15; 34B18, 35J25,
35J61.}  \thanks{Partially supported by Secyt-UNC and CONICET. The first
author wants to dedicate this work to his teacher and friend Tom\'{a}s Godoy.} }
\author{U. Kaufmann, I. Medri\thanks{\textit{E-mail addresses. }%
kaufmann@mate.uncor.edu (U. Kaufmann, Corresponding Author),
medri@mate.uncor.edu (I. Medri).}
\and \noindent\\{\small FaMAF, Universidad Nacional de C\'{o}rdoba, (5000) C\'{o}rdoba,
Argentina}}
\maketitle

\begin{abstract}
We study existence and nonexistence of strictly positive solutions for the
elliptic problems of the form $Lu=m\left(  x\right)  u^{p}$ in a bounded open
interval, with zero boundary conditions, where $L$ is a strongly uniformly
elliptic differential operator, $p\in\left(  0,1\right)  $, and $m$ is a
function that changes sign. We also characterize the set of values $p$ for
which the problem admits a solution, and in addition an existence result for
other nonlinearities is presented.

\end{abstract}

\section{Introduction}

For $\alpha<\beta$, let $\Omega:=(\alpha,\beta)$ and $m\in L^{2}\left(
\Omega\right)  $ be a function that changes sign in $\Omega$. Let $p\in\left(
0,1\right)  $ and let $L$ be a one-dimensional strongly uniformly elliptic
differential operator given by
\begin{equation}
Lu:=-a(x)u^{\prime\prime}+b(x)u^{\prime}+c(x)u, \label{L}%
\end{equation}
where $a,b\in C\left(  \overline{\Omega}\right)  $, $0\leq c\in L^{\infty
}(\Omega)$ and $a\left(  x\right)  \geq\lambda>0$ for all $x\in\Omega$. Our
aim in the present paper is to consider the matter of existence and
nonexistence of solutions for problems of the form%

\begin{equation}
\left\{
\begin{array}
[c]{ll}%
Lu=mu^{p} & \text{in }\Omega\\
u>0 & \text{in }\Omega\\
u=0 & \text{on }\partial\Omega.
\end{array}
\right.  \label{prob}%
\end{equation}

The question of existence of strictly positive solutions for semilinear
Dirichlet problems with indefinite nonlinearities as (\ref{prob}) is
challenging and intriguing, and to our knowledge there are few results
concerning this issue. In contrast to superlinear problems where any
nonnegative (and nontrivial) solution is automatically positive (and in fact
is in the interior of the positive cone under standard assumptions), for the
analogous sublinear equations the situation is far less clear, even in the
one-dimensional case. For instance, it is known that if $m$ is smooth and
$m^{+}\not \equiv 0$ then for any $p\in\left(  0,1\right)  $ there exist
nontrivial nonnegative solutions that actually vanish in a subset of $\Omega$
(see e.g. \cite{bandle}, \cite{publi}), and when $L=-u^{\prime\prime}$ one may
also construct examples of strictly positive solutions that do not belong to
the interior of the positive cone (see \cite{ultimo}).

The problem (\ref{prob}) was considered recently in \cite{ultimo} for the
laplacian operator, where several \textit{non-comparable }sufficient
conditions for the existence of solutions where proved under some evenness
assumptions on $m$. In the present paper we shall adapt and extend the
approach in \cite{ultimo} in order to derive our main results for a general
operator. More precisely, in Section 3 we shall give two non-comparable
sufficient conditions on $m$ in the case $b\equiv0$ (see Theorem \ref{aa} and
Remark \ref{sepo}), and when $b\not \equiv 0$ we shall also exhibit sufficient
conditions in Theorem \ref{bien} and Corollary \ref{puf}. Let us mention that
these last conditions are non-comparable between each other nor between the
ones in Theorem \ref{aa}. Moreover, one of them substantially improves the
results known for $L=-u^{\prime\prime}$ (see Remarks \ref{ojito1},
\ref{ojito2} and \ref{lm}). Also, as a consequence of the aforementioned
results we shall characterize the set of $p^{\prime}$s such that (\ref{prob})
admits a solution and we shall deduce an existence theorem for other
nonlinearities (see Corollaries \ref{ppp} and \ref{fff} respectively). Let us
finally say that necessary conditions on $m$ for the existence of solutions
are stated in Theorem \ref{necee}.

In order to relate our results to others already existing let us mention that
to our knowledge no necessary condition on $m$ is known in the case of a
general operator (other than the obvious one derived from the maximum
principle, i.e. $m^{+}\not \equiv 0$), and the only sufficient condition we
found in the literature is that the solution $\varphi$ of $L\varphi=m$ in
$\Omega$, $\varphi=0$ on $\partial\Omega$, satisfies $\varphi>0$ in $\Omega$
(see \cite{jesusultimo}, Theorem 4.4, or \cite{hand}, Theorem 10.6).\ Let us
note that although the above condition is even true for the $n$-dimensional
problem, it is far from being necessary in the sense that there are examples
of (\ref{prob}) having a solution but with the corresponding $\varphi$
satisfying $\varphi<0$ in $\Omega$ (cf. \cite{ultimo}). Concerning the
laplacian operator, (\ref{prob}) was treated in Theorem 2.1 in \cite{ultimo},
and as we said before there are also further results there under different
evenness assumptions on $m$. Let us finally mention that existence of
solutions for problem (\ref{prob}) has also been studied when $L=-u^{\prime
\prime}$ and $m\geq0$ but assuming that $m\in C\left(  \Omega\right)  $ (see
e.g. \cite{zhang}, \cite{chapa} and the references therein), and also some
similar results to the ones that appear here have been obtained recently by
the authors in \cite{aust} for some related problems involving quasilinear operators.

We would like to conclude this introduction with some few words on the
corresponding $n$-dimensional problem. As we noticed in the above paragraph
the condition in \cite{jesusultimo} is still valid in this case, and some of
the techniques in \cite{ultimo} can be applied if $L=-\Delta$ (see Section 3
in \cite{ultimo} for the radial case, and also \cite{junows}). We are strongly
convinced that some of the theorems presented here should still have some
counterpart in $n$ dimensions but we are not able to provide a proof.

\section{Preliminaries and auxiliary results}

Since $a\left(  x\right)  \geq\lambda>0$ for all $x\in\Omega$ and $a\in
C\left(  \overline{\Omega}\right)  $, from now on we consider without loss of
generality that $L$ is given by%
\begin{equation}
Lu:=-u^{\prime\prime}+b(x)u^{\prime}+c(x)u\text{,} \label{L2}%
\end{equation}
with $b$ and $c$ as in (\ref{L}). For $f\in L^{r}(\Omega)$ with $r>1$ we say
that $u$ is a (strong) solution of the problem $Lu=f$ in $\Omega$, $u=0$ in
$\partial\Omega$, if $u\in W^{2,r}(\Omega)\cap W_{0}^{1,r}\left(
\Omega\right)  $ and the equation is satisfied $a.e.$ $x\in\Omega$. Given
$g:\Omega\times\mathbb{R}\rightarrow\mathbb{R}$ a Carathe\'{o}dory function
such that $g\left(  .,\xi\right)  \in L^{2}(\Omega)$ for all $\xi$, we say
that $u$ is a (weak) subsolution of
\begin{equation}
\left\{
\begin{array}
[c]{ll}%
Lu=g\left(  x,u\right)  & \text{in }\Omega\\
u=0 & \text{on }\partial\Omega
\end{array}
\right.  \label{g}%
\end{equation}
if $u\in W^{1,2}\left(  \Omega\right)  $, $u\leq0$ on $\partial\Omega$ and
\[
\int_{\Omega}u^{\prime}\phi^{\prime}+bu^{\prime}\phi+cu\phi\leq\int_{\Omega
}g\left(  x,u\right)  \phi\qquad\text{for all }0\leq\phi\in W_{0}^{1,2}%
(\Omega)\text{.}%
\]
(Weak) supersolutions are defined analogously.

The following lemma is a direct consequence of the integration by parts
formula (e.g. \cite{brezislibro}, Corollary 8.10).

\begin{lemma}
\label{coro}For $i:1,...,n$, let $u_{i}\in W^{2,2}(x_{i},x_{i+1})$ or
$u_{i}\in C^{2}(x_{i},x_{i+1})\cap C^{1}\left(  \left[  x_{i},x_{i+1}\right]
\right)  $ such that $u_{i}(x_{i+1})=u_{i+1}(x_{i+1})$, $u_{i}^{\prime
}(x_{i+1})\leq u_{i+1}^{\prime}(x_{i+1})$ and
\[
-u_{i}^{\prime\prime}+bu_{i}^{\prime}+cu_{i}\leq g\left(  x,u_{i}\right)
\qquad a.e.\text{ }x\in(x_{i},x_{i+1})\text{ for all\ }i:1,...,n.
\]
Let $\Omega:=\left(  x_{1},x_{n+1}\right)  $ and set $u(x):=u_{i}(x)$ for all
$x\in\Omega$. Then $u\in W^{1,2}(\Omega)$ and%
\[
\int_{\Omega}u^{\prime}\phi^{\prime}+bu^{\prime}\phi+cu\phi\leq\int_{\Omega
}g\left(  x,u\right)  \phi\text{\qquad for all }0\leq\phi\in W_{0}%
^{1,2}(\Omega)\text{.}%
\]
In particular, if also $u\leq0$ on $\partial\Omega$, then $u$ is a subsolution
of (\ref{g}).
\end{lemma}

The next remark compiles some necessary facts about problem (\ref{prob}).

\begin{remark}
\label{homsup} (i) It is immediate to check that (\ref{prob}) possesses a
solution if and only if it has a solution with $\tau m$ in place of $m$, for
any $\tau>0$. \newline(ii) Let us write as usual $m=m^{+}-m^{-}$ with
$m^{+}=\max\left(  m,0\right)  $ and $m^{-}=\max\left(  -m,0\right)  $. It is
also easy to verify that (\ref{prob}) admits arbitrarily large supersolutions
(if $m^{+}\not \equiv 0$; if $m^{+}\equiv0$ there is no solution by the
maximum principle). Indeed, let $\varphi>0$ be the solution of $L\varphi
=m^{+}$ in $\Omega$, $\varphi=0$ on $\partial\Omega$. Let $k\geq(\left\Vert
\varphi\right\Vert _{\infty}+1)^{p/(1-p)}$. Then $k(\varphi+1)$ is a
supersolution since
\begin{equation}
L(k(\varphi+1))\geq kL\varphi\geq(k(\left\Vert \varphi\right\Vert _{\infty
}+1))^{p}m^{+}\geq(k(\varphi+1))^{p}m\text{\qquad in }\Omega\label{ju}%
\end{equation}
and $\varphi=k>0$ on $\partial\Omega$. $\blacksquare$
\end{remark}

The two following lemmas provide some useful upper bounds for the $L^{\infty}%
$-norm of the nonnegative subsolutions of (\ref{prob}). In order to avoid
overloading the notation we write from now on
\[
\overline{B}_{\alpha}\left(  x\right)  :=e^{\int_{\alpha}^{x}b\left(
r\right)  dr},\qquad\underline{B}_{\alpha}\left(  x\right)  :=e^{-\int
_{\alpha}^{x}b\left(  r\right)  dr}.
\]

\begin{lemma}
\label{inerte}Let $0\leq u\in W^{2,2}\left(  \Omega\right)  $ be such that
$Lu\leq mu^{p}$ in $\Omega$. Then%
\begin{equation}
\left\Vert u\right\Vert _{L^{\infty}(\Omega)}\leq\left[  \int_{\alpha}^{\beta
}\overline{B}_{\alpha}\left(  x\right)  \left\Vert m^{+}\underline{B}_{\alpha
}\right\Vert _{L^{1}(\alpha,x)}dx\right]  ^{1/\left(  1-p\right)  }.
\label{iner}%
\end{equation}

\end{lemma}

\textit{Proof}. Since $\underline{B}_{\alpha},u^{\prime}\in W^{1,2}(\Omega)$,
we may apply the product differentiation rule and hence%
\begin{gather*}
-(\underline{B}_{\alpha}u^{\prime})^{\prime}\leq-(\underline{B}_{\alpha
}u^{\prime})^{\prime}+\underline{B}_{\alpha}cu=\underline{B}_{\alpha
}(-u^{\prime\prime}+bu^{\prime}+cu)\leq\\
\underline{B}_{\alpha}mu^{p}\leq\underline{B}_{\alpha}m^{+}\left\Vert
u\right\Vert _{L^{\infty}(\Omega)}^{p}.
\end{gather*}
Integrating on $\left(  \alpha,x\right)  $ for $x\in(\alpha,\beta)$ (see e.g.
\cite{brezislibro}, Theorem 8.2) and noting that $\underline{B}_{\alpha
}(\alpha)u^{\prime}\left(  \alpha\right)  =u^{\prime}\left(  \alpha\right)
\geq0$ we obtain
\[
-\underline{B}_{\alpha}(x)u^{\prime}(x)\leq\left\Vert u\right\Vert
_{L^{\infty}(\Omega)}^{p}\int_{\alpha}^{x}\underline{B}_{\alpha}%
(t)m^{+}(t)dt.
\]
Dividing by $\underline{B}_{\alpha}(x)>0$ and integrating now on $\left(
y,\beta\right)  $ for $y\in(\alpha,\beta)$, since $u(\beta)=0$ we get
\[
0\leq\frac{u(y)}{\left\Vert u\right\Vert _{L^{\infty}(\Omega)}^{p}}\leq
\int_{y}^{\beta}\left[  \overline{B}_{\alpha}(x)\int_{\alpha}^{x}\underline
{B}_{\alpha}(t)m^{+}(t)dt\right]  dx\qquad\text{for all }y\in(\alpha,\beta),
\]
and the lemma follows. $\blacksquare$

\qquad

Let
\begin{equation}
M^{+}:=\left\{  x\in\Omega:m\geq0\right\}  ,\qquad M^{-}:=\left\{  x\in
\Omega:m<0\right\}  . \label{MM}%
\end{equation}

\begin{lemma}
\label{dudu}Let $0\leq u\in W^{2,2}\left(  \Omega\right)  $ be such that
$Lu\leq mu^{p}$ in $\Omega$, and let $M^{+}$ be given by (\ref{MM}). If $c>0$
in $M^{+}$, then%
\[
\left\Vert u\right\Vert _{L^{\infty}(\Omega)}\leq\left[  \sup_{x\in M^{+}%
}\frac{m^{+}\left(  x\right)  }{c\left(  x\right)  }\right]  ^{1/\left(
1-p\right)  }.
\]

\end{lemma}

\textit{Proof. }Without loss of generality we assume that $u\not \equiv 0$.
Furthermore, let us suppose first that $\left\Vert u\right\Vert _{L^{\infty
}(\Omega)}>1$. Let $x_{0}\in\Omega$ be a point where $u$ attains its absolute
maximum. There exists $\delta>0$ such that $u\geq1$ in $I_{\delta}\left(
x_{0}\right)  :=\left(  x_{0}-\delta,x_{0}+\delta\right)  $. There also exist
$x_{1},x_{2}\in I_{\delta}\left(  x_{0}\right)  $ satisfying $x_{1}%
<x_{0}<x_{2}$ and $u^{\prime}\left(  x_{2}\right)  \leq0\leq u^{\prime}\left(
x_{1}\right)  $. We have that
\[
-\left(  \underline{B}_{\alpha}u^{\prime}\right)  ^{\prime}+\underline
{B}_{\alpha}cu\leq\underline{B}_{\alpha}mu^{p}\leq\underline{B}_{\alpha}%
m^{+}u^{p}\qquad\text{in }\Omega
\]
and so in $I_{\delta}\left(  x_{0}\right)  $ we get that (because $u\geq1$ in
$I_{\delta}\left(  x_{0}\right)  $) $-\left(  \underline{B}_{\alpha}u^{\prime
}\right)  ^{\prime}\leq\underline{B}_{\alpha}\left(  m^{+}-c\right)  u$.
Integrating on $\left(  x_{1},x_{2}\right)  $ we derive%
\begin{equation}
0\leq\underline{B}_{\alpha}\left(  x_{1}\right)  u^{\prime}\left(
x_{1}\right)  -\underline{B}_{\alpha}\left(  x_{2}\right)  u^{\prime}\left(
x_{2}\right)  =\int_{x_{1}}^{x_{2}}-\left(  \underline{B}_{\alpha}u^{\prime
}\right)  ^{\prime}\leq\int_{x_{1}}^{x_{2}}\underline{B}_{\alpha}\left(
m^{+}-c\right)  u. \label{cero}%
\end{equation}
Since $u\geq1$ in $\left(  x_{1},x_{2}\right)  $ and $\underline{B}_{\alpha
}\geq e^{-\left\Vert b^{+}\right\Vert _{\infty}\left(  x_{2}-\alpha\right)  }$
in $\left(  x_{1},x_{2}\right)  $, from (\ref{cero}) it follows that there
exists $E\subset\left(  x_{1},x_{2}\right)  $ with $\left\vert E\right\vert
>0$ (where $\left\vert E\right\vert $ denotes the Lebesgue measure of $E$)
such that $m^{+}\left(  x\right)  \geq c\left(  x\right)  $ $a.e.$ $x\in E$.
Moreover, due to the fact that $c>0$ $a.e.$ $x\in M^{+}$ it must hold that
$m^{+}>0$ $a.e.$ $x\in E$. In particular, $E\subset M^{+}$ and therefore%
\begin{equation}
1\leq\sup_{x\in E}\frac{m^{+}\left(  x\right)  }{c\left(  x\right)  }\leq
\sup_{x\in M^{+}}\frac{m^{+}\left(  x\right)  }{c\left(  x\right)  }.
\label{casi}%
\end{equation}

Let $u$ now be as in the statement of the lemma, and let $\varepsilon>0$. Then%
\[
L\frac{u}{\left\Vert u\right\Vert _{\infty}-\varepsilon}\leq\frac{m}{\left(
\left\Vert u\right\Vert _{\infty}-\varepsilon\right)  ^{1-p}}\left(  \frac
{u}{\left\Vert u\right\Vert _{\infty}-\varepsilon}\right)  ^{p}.
\]
Applying the first part of the proof with $m/\left(  \left\Vert u\right\Vert
_{\infty}-\varepsilon\right)  ^{1-p}$ and $u/\left(  \left\Vert u\right\Vert
_{\infty}-\varepsilon\right)  $ in place of $m$ and $u$ respectively, from
(\ref{casi}) we deduce that
\[
\left(  \left\Vert u\right\Vert _{L^{\infty}(\Omega)}-\varepsilon\right)
^{1-p}\leq\sup_{x\in M^{+}}\frac{m^{+}\left(  x\right)  }{c\left(  x\right)  }%
\]
and since $\varepsilon$ is arbitrary this ends the proof of the lemma.
$\blacksquare$

\qquad

We shall need the next result when we characterize the set of $p^{\prime}$s
such that (\ref{prob}) admits a solution.

\begin{lemma}
\label{qp} Suppose (\ref{prob}) has a solution $u\in W^{2,2}\left(
\Omega\right)  $, and let $q\in\left(  p,1\right)  $. Then there exists $v\in
W^{2,2}\left(  \Omega\right)  $ solution of (\ref{prob}) with $q$ in place of
$p$.
\end{lemma}

\textit{Proof}. Let $\gamma:=\left(  1-p\right)  /\left(  1-q\right)  $. Let
$0\leq\phi\in C_{c}^{\infty}\left(  \Omega\right)  $, and let $\Omega^{\prime
}$ be an open set such that \textit{supp} $\phi\subset\Omega^{\prime}%
\Subset\Omega$. One can check that $u^{\gamma}\in W_{0}^{1,2}\left(
\Omega\right)  \cap W^{2,2}\left(  \Omega^{\prime}\right)  $. Furthermore,
noticing that $\gamma>1$ and $\gamma-1+p=\gamma q$ we find that
\begin{gather*}
L\left(  u^{\gamma}\right)  =-\gamma\left(  u^{\prime\prime}u^{\gamma
-1}+\left(  \gamma-1\right)  u^{\gamma-2}\left(  u^{\prime}\right)
^{2}\right)  +b\gamma u^{\gamma-1}u^{\prime}+cu^{\gamma}\leq\\
\gamma u^{\gamma-1}\left(  -u^{\prime\prime}+bu^{\prime}+cu\right)  \leq\gamma
u^{\gamma-1}mu^{p}=\gamma m\left(  u^{\gamma}\right)  ^{q}\qquad\text{in
}\Omega^{\prime}\text{.}%
\end{gather*}
Multiplying the above inequality by $\phi$, integrating over $\Omega^{\prime}$
and using the integration by parts formula we obtain that%
\begin{gather*}
\int_{\Omega}\left(  u^{\gamma}\right)  ^{\prime}\phi^{\prime}+b\left(
u^{\gamma}\right)  ^{\prime}\phi+cu^{\gamma}\phi=\int_{\Omega^{\prime}}\left[
-\left(  u^{\gamma}\right)  ^{\prime\prime}+b\left(  u^{\gamma}\right)
^{\prime}+cu^{\gamma}\right]  \phi\leq\\
\gamma\int_{\Omega}m\left(  u^{\gamma}\right)  ^{q}\phi.
\end{gather*}
Now, let $0\leq v\in W_{0}^{1,2}\left(  \Omega\right)  $. There exists
$\left\{  \phi_{n}\right\}  _{n\in\mathbb{N}}\subset C_{c}^{\infty}\left(
\Omega\right)  $ with $\phi_{n}\geq0$ in $\Omega$ and such that $\phi
_{n}\rightarrow v$ in $W^{1,2}\left(  \Omega\right)  $ (e.g. \cite{chipot}, p.
50). Employing the above inequality with $\phi_{n}$ in place of $\phi$ and
going to the limit we see that $u^{\gamma}$ is a subsolution of (\ref{prob})
with $\gamma m$ in place of $m$. Thus, taking into account Remark \ref{homsup}
(i) and (ii) we get a solution $v\in W_{0}^{1,2}\left(  \Omega\right)  $ of
(\ref{prob}), and by standard regularity arguments $v\in W^{2,2}\left(
\Omega\right)  $. $\blacksquare$

\section{Main results}

We set
\begin{equation}
C_{p}:=\frac{2\left(  1+p\right)  }{\left(  1-p\right)  ^{2}},\label{cp}%
\end{equation}
and for any interval $I$,%
\[
\lambda_{1}\left(  m,I\right)  :=\text{the positive principal eigenvalue for
}m\text{ in }I\text{.}%
\]

\begin{theorem}
\label{aa}Assume $b\equiv0$. Let $m\in L^{2}(\Omega)$ with $m^{-}\in
L^{\infty}(\Omega)$ and suppose there exist $\alpha\leq x_{0}<x_{1}\leq\beta$
such that $0\not \equiv m\geq0$ in $I:=(x_{0},x_{1})$. Let $\gamma
:=\max\left\{  (\beta-x_{0}),(x_{1}-\alpha)\right\}  $ and let $C_{p}$ be
given by (\ref{cp}).\newline(i) If it holds that
\begin{equation}
\frac{\left\Vert m^{-}\right\Vert _{L^{\infty}(\Omega)}}{\left\Vert
c\right\Vert _{L^{\infty}(\Omega)}}\sinh^{2}\left[  \gamma\sqrt{\frac
{\left\Vert c\right\Vert _{\infty}}{C_{p}}}\right]  \leq\frac{1}{\lambda
_{1}(m,I)} \label{seno}%
\end{equation}
then there \textit{exists }$u\in W^{2,2}(\Omega)$\textit{\ solution of
}(\ref{prob})\textit{.}\newline(ii) If it holds that
\begin{equation}
\frac{\left\Vert m^{-}\right\Vert _{L^{\infty}(\Omega)}}{\left\Vert
c\right\Vert _{L^{\infty}(\Omega)}}\left[  \cosh\left(  \gamma\sqrt{\left(
1-p\right)  \left\Vert c\right\Vert _{L^{\infty}(\Omega)}}\right)  -1\right]
\leq\frac{1}{\lambda_{1}(m,I)} \label{expo}%
\end{equation}
then there \textit{exists }$u\in W^{2,2}(\Omega)$\textit{\ solution of
}(\ref{prob})\textit{.}\newline
\end{theorem}

\textit{Proof. }Recalling Remark \ref{homsup} it suffices to construct a
strictly positive (in $\Omega$) subsolution $u$ for (\ref{prob}) with $\tau m$
in place of $m$, for some $\tau>0$. Moreover, without loss of generality we
may assume that $\alpha<x_{0}<x_{1}<\beta$ (in fact, it shall be clear from
the proof how to proceed if either $x_{0}=\alpha$ or $x_{1}=\beta$). In order
to provide such $u$ we shall employ Lemma \ref{coro} with $n=3$ and $g\left(
x,\xi\right)  =\tau m\left(  x\right)  \xi^{p}$.

We shall take $u_{2}>0$ with $\left\Vert u_{2}\right\Vert _{L^{\infty}\left(
I\right)  }=1$ as the positive principal eigenfunction associated to the
weight $m$ in $I$, that is satisfying%
\[
\left\{
\begin{array}
[c]{ll}%
Lu_{2}=\lambda_{1}(m,I)mu_{2} & \text{in }I\\
u_{2}=0 & \text{on }\partial I.
\end{array}
\right.
\]
Since $m\geq0$ in $I$, for $\tau>0$ we have that $Lu_{2}=\lambda
_{1}(m,I)mu_{2}\leq\tau mu_{2}^{p}$ whenever
\begin{equation}
\lambda_{1}(m,I)\leq\tau.\label{tri}%
\end{equation}
On the other hand, suppose now that (\ref{seno}) holds and pick $\tau$
satisfying
\begin{equation}
\frac{\left\Vert m^{-}\right\Vert _{L^{\infty}(\Omega)}}{\left\Vert
c\right\Vert _{L^{\infty}(\Omega)}}\sinh^{2}\left[  \gamma\sqrt{\frac
{\left\Vert c\right\Vert _{\infty}}{C_{p}}}\right]  \leq\frac{1}{\tau}%
\leq\frac{1}{\lambda_{1}(m,I)}\label{tau}%
\end{equation}
(in particular, (\ref{tri}) holds). Let $x\in\left[  \alpha,x_{1}\right]  $
and define
\[
f(x)=\sqrt{\frac{\tau\left\Vert m^{-}\right\Vert _{\infty}}{\left\Vert
c\right\Vert _{\infty}}}\sinh\left[  \sqrt{\frac{\left\Vert c\right\Vert
_{\infty}}{C_{p}}}\left(  x-\alpha\right)  \right]  .
\]
A few computations show that $C_{p}\left(  f^{\prime}\right)  ^{2}-\left\Vert
c\right\Vert _{\infty}f^{2}=\tau\left\Vert m^{-}\right\Vert _{\infty}$ in
$(\alpha,x_{1})$. Moreover, $f(\alpha)=0$, $f\left(  x\right)  >0$ for
$x\in\left(  \alpha,x_{1}\right)  $ and $f^{\prime},f^{\prime\prime}\geq0$ for
such $x$. Let us now fix $k:=2/\left(  1-p\right)  $. Then we have
\begin{equation}
kp=k-2\text{,}\qquad k\left(  k-1\right)  =C_{p}.\label{kk}%
\end{equation}
We set $u_{1}:=f^{k}$. Taking into account (\ref{kk}) and the above mentioned
facts we find that
\begin{gather}
Lu_{1}=-k\left[  \left(  k-1\right)  f^{k-2}\left(  f^{\prime}\right)
^{2}+f^{k-1}f^{\prime\prime}\right]  +cf^{k}\leq\label{ups}\\
-C_{p}f^{k-2}\left(  f^{\prime}\right)  ^{2}+\left\Vert c\right\Vert _{\infty
}f^{k}=-f^{k-2}\tau\left\Vert m^{-}\right\Vert _{\infty}\leq\tau mu_{1}%
^{p}\qquad\text{in }\left(  \alpha,x_{1}\right)  .\nonumber
\end{gather}
Furthermore, since $f$ is increasing we get that $\left\Vert u_{1}\right\Vert
_{\infty}=\left[  f(x_{1})\right]  ^{k}$ and therefore using the first
inequality in (\ref{tau}) and the fact that $x_{1}-\alpha\leq\gamma$ one can
verify that $\left\Vert u_{1}\right\Vert _{\infty}\leq1$.

In a similar way, if for $x\in\left[  x_{0},\beta\right]  $ we define
$u_{3}:=g^{k}$ where $g$ is given by%
\[
g(x):=\sqrt{\frac{\tau\left\Vert m^{-}\right\Vert _{\infty}}{\left\Vert
c\right\Vert _{\infty}}}\sinh\left[  \sqrt{\frac{\left\Vert c\right\Vert
_{\infty}}{C_{p}}}\left(  \beta-x\right)  \right]  ,
\]
then $Lu_{3}\leq\tau mu_{3}^{p}$ in $\left(  x_{0},\beta\right)  $,
$\left\Vert u_{3}\right\Vert _{\infty}\leq1$, $u_{3}(\beta)=0$ and
$u_{3}\left(  x\right)  >0$ for $x\in\left(  x_{0},\beta\right)  $.

We choose now
\begin{gather*}
\underline{x}_{0}:=\sup\left\{  x\in I:u_{1}\left(  y\right)  >u_{2}\left(
y\right)  \text{ for all }y\in\left(  x_{0},x\right]  \right\}  ,\\
\overline{y}:=\max\left\{  x\in I:u_{2}\left(  x\right)  =1\right\}
,\qquad\underline{y}:=\min\left\{  x\in I:u_{2}\left(  x\right)  =1\right\}  .
\end{gather*}
We observe that $\underline{x}_{0}\in I$ exists because $u_{1}\left(
\alpha\right)  =u_{2}\left(  x_{0}\right)  =0$ and $u_{1}(x_{1})\leq
1=\left\Vert u_{2}\right\Vert _{\infty}$. Moreover, since $u_{1}$ and $u_{2}$
are $C^{1}$, by the definition of $\underline{x}_{0}$ we have that
$u_{1}(\underline{x}_{0})=u_{2}(\underline{x}_{0})$ and $u_{1}^{\prime
}(\underline{x}_{0})\leq u_{2}^{\prime}(\underline{x}_{0})$ (for the last
inequality it is enough to note that $\frac{u_{1}(x)-u_{1}(\underline{x}_{0}%
)}{x-\underline{x}_{0}}<\frac{u_{2}(x)-u_{2}(\underline{x}_{0})}%
{x-\underline{x}_{0}}$ for every $x\in\left(  x_{0},\underline{x}_{0}\right)
$), and also clearly $\underline{x}_{0}<\underline{y}$. Analogously, there
exists $\overline{x}_{1}\in I$ such that $u_{2}\left(  \overline{x}%
_{1}\right)  =u_{3}\left(  \overline{x}_{1}\right)  $ and $u_{2}^{\prime
}(\overline{x}_{1})\leq u_{3}^{\prime}(\overline{x}_{1})$, and satisfying
$\overline{x}_{1}>\overline{y}$. In particular, $\underline{x}_{0}%
<\overline{x}_{1}$. Hence, defining $u$ by $u:=u_{1}$ in $\left[
\alpha,\underline{x}_{0}\right]  $, $u:=u_{2}$ in $\left[  \underline{x}%
_{0},\overline{x}_{1}\right]  $ and $u:=u_{3}$ in $\left[  \overline{x}%
_{1},\beta\right]  $, we have that $u=0$ on $\partial\Omega$ and $u$ fulfills
the hypothesis of Lemma \ref{coro} and as we said before this proves (i) (let
us mention that if $x_{0}=\alpha$ then in order to build $u$ we only use
$u_{2}$ and $u_{3}$, and if $x_{1}=\beta$ then we do not need $u_{3}$).

Let us prove (ii). We shall take $u_{2}$ as above. We now fix $\tau$ such
that
\begin{equation}
\frac{\left\Vert m^{-}\right\Vert _{L^{\infty}(\Omega)}}{\left\Vert
c\right\Vert _{L^{\infty}(\Omega)}}\left[  \cosh\left(  \gamma\sqrt{\left(
1-p\right)  \left\Vert c\right\Vert _{L^{\infty}(\Omega)}}\right)  -1\right]
\leq\frac{1}{\tau}\leq\frac{1}{\lambda_{1}(m,I)}.\label{tauu}%
\end{equation}
We set $k:=1/\left(  1-p\right)  $, and for $x\in\left[  \alpha,x_{1}\right]
$ we define%
\[
f\left(  x\right)  :=\frac{\tau\left\Vert m^{-}\right\Vert _{\infty}%
}{\left\Vert c\right\Vert _{\infty}}\left[  \cosh\left(  \sqrt{\frac
{\left\Vert c\right\Vert _{\infty}}{k}}(x-\alpha)\right)  -1\right]  .
\]
Then $f(\alpha)=0$, $f>0$ in $(\alpha,x_{1})$ and $f^{\prime}\geq0$.
Furthermore, by the first inequality in (\ref{tauu}) $\left\Vert
u_{1}\right\Vert _{\infty}\leq1$, and it can be seen that $kf^{\prime\prime
}-\left\Vert c\right\Vert _{\infty}f=\tau\left\Vert m^{-}\right\Vert _{\infty
}$. Define now $u_{1}:=f^{k}$. Observing that $kp=k-1$ we derive that
\begin{gather*}
Lu_{1}=-k\left[  \left(  k-1\right)  f^{k-2}\left(  f^{\prime}\right)
^{2}+f^{k-1}f^{\prime\prime}\right]  +cf^{k}\leq\\
-kf^{k-1}f^{\prime\prime}+\left\Vert c\right\Vert _{\infty}f^{k}=-f^{k-1}%
\tau\left\Vert m^{-}\right\Vert _{\infty}\leq\tau mu_{1}^{p}\qquad\text{in
}\left(  \alpha,x_{1}\right)  .
\end{gather*}
In the same way, if for $x\in\left[  x_{0},\beta\right]  $ we set
$u_{3}:=g^{k}$ where $g$ is given by%
\[
g(x):=\frac{\tau\left\Vert m^{-}\right\Vert _{\infty}}{\left\Vert c\right\Vert
_{\infty}}\left[  \cosh\left(  \sqrt{\frac{\left\Vert c\right\Vert _{\infty}%
}{k}}(\beta-x)\right)  -1\right]  ,
\]
then $Lu_{3}\leq\tau mu_{3}^{p}$ in $\left(  x_{0},\beta\right)  $,
$\left\Vert u_{3}\right\Vert _{\infty}\leq1$, $u_{3}(\beta)=0$ and $u_{3}>0$
in $\left(  x_{0},\beta\right)  $. Now the proof of (ii) can be finished as in
(i). $\blacksquare$

\begin{remark}
\label{sepo}Let us mention that the inequalities in (i) and (ii) are not
comparable. Indeed, we first check that for $p\approx1$ (\ref{seno}) is better
than (\ref{expo}). Let $\kappa:=\gamma\sqrt{\left\Vert c\right\Vert _{\infty}%
}$. Since $\frac{1}{\sqrt{C_{p}}}=\left(  1-p\right)  \sqrt{\frac{1}{2\left(
1+p\right)  }}$, it is enough to observe that
\[
0\leq\lim_{p\rightarrow1^{-}}\frac{\sinh^{2}\left[  \kappa\left(  1-p\right)
\sqrt{\frac{1}{2\left(  1+p\right)  }}\right]  }{\cosh\left(  \kappa\sqrt
{1-p}\right)  -1}\leq\lim_{p\rightarrow1^{-}}\frac{\sinh^{2}\left(
\kappa\left(  1-p\right)  \right)  }{\cosh\left(  \kappa\sqrt{1-p}\right)
-1}=0.
\]
We now show that for $0<p\approx0$ (\ref{expo}) is better than (\ref{seno}).
It suffices to prove this for $p=0$ because the dependence on $p$ in both
inequalities is continuous. For $p=0$ (\ref{seno}) and (\ref{expo}) become%
\begin{gather*}
\frac{\left\Vert m^{-}\right\Vert _{\infty}}{\left\Vert c\right\Vert _{\infty
}}\sinh^{2}\left(  \kappa/\sqrt{2}\right)  \leq\frac{1}{\lambda_{1}(m,I)}\\
\frac{\left\Vert m^{-}\right\Vert _{\infty}}{\left\Vert c\right\Vert _{\infty
}}\left(  \cosh\kappa-1\right)  \leq\frac{1}{\lambda_{1}(m,I)}%
\end{gather*}
and so we only have to check that for every $x>0$ it holds that $\sinh
^{2}\left(  x/\sqrt{2}\right)  >\cosh x-1$ which is easy to verify.
$\blacksquare$
\end{remark}

\begin{remark}
\label{lapla}If in (\ref{seno}) we take limit as $\left\Vert c\right\Vert
_{L^{\infty}(\Omega)}\rightarrow0$ we arrive to the condition
\begin{equation}
\frac{\gamma^{2}}{C_{p}}\left\Vert m^{-}\right\Vert _{L^{\infty}(\Omega)}%
\leq\frac{1}{\lambda_{1}(m,I)} \label{lap}%
\end{equation}
which is the one that appears for $L=-u^{\prime\prime}$ in Theorem 2.1 in
\cite{ultimo}. $\blacksquare$
\end{remark}

\begin{remark}
In the statement of Theorem \ref{aa} one can replace the condition
(\ref{seno}) by%
\begin{gather}
\frac{\left\Vert m^{-}\right\Vert _{L^{\infty}(\Omega)}}{\left\Vert
c\right\Vert _{L^{\infty}(M^{-})}}\sinh^{2}\left[  \gamma\sqrt{\frac
{\left\Vert c\right\Vert _{L^{\infty}(M^{-})}}{C_{p}}}\right]  \leq\frac
{1}{\lambda_{1}(m,I)}\qquad\text{and}\label{rem}\\
c\leq m^{+}\text{ in }M^{+}, \label{remm}%
\end{gather}
where $M^{+}$ and $M^{-}$ are given by (\ref{MM}). Indeed, we first observe
that if (\ref{rem}) holds then one can reason as in (\ref{ups}) and prove that
$Lu_{1}\leq\tau mu_{1}^{p}$ in $\left(  x_{0},\beta\right)  \cap M^{-}$. On
the other side, if (\ref{remm}) is true then since in the proof of the theorem
$f$ is chosen satisfying $f^{\prime\prime}\geq0$ and $\left\Vert
f^{k}\right\Vert _{\infty}\leq1$, then we also have
\begin{gather*}
Lu_{1}=-k\left[  \left(  k-1\right)  f^{k-2}\left(  f^{\prime}\right)
^{2}+f^{k-1}f^{\prime\prime}\right]  +cf^{k}\leq\\
cf^{k}\leq m^{+}f^{k}\leq m^{+}f^{kp}=mu_{1}^{p}\qquad\text{in }\left(
x_{0},\beta\right)  \cap M^{+}.
\end{gather*}
The same reasoning can be done for $u_{3}$ and hence the proof can be
continued as in the theorem. A similar observation is valid for (\ref{expo}).
$\blacksquare$
\end{remark}

\begin{theorem}
\label{bien}Let $m\in L^{2}(\Omega)$ and suppose there exist $\alpha\leq
x_{0}<x_{1}\leq\beta$ such that $0\not \equiv m\geq0$ in $I:=(x_{0},x_{1})$.
Let $C_{p}$ be given by (\ref{cp}).\newline(i) If $m^{-}\in L^{\infty}%
(\Omega)$ and it holds that%
\begin{gather}
0<\frac{\left(  \gamma_{b}\left\Vert \underline{B}_{\alpha}\right\Vert
_{L^{\infty}(\Omega)}\right)  ^{2}}{C_{p}-\left\Vert c\right\Vert _{L^{\infty
}(\Omega)}\left(  \gamma_{b}\left\Vert \underline{B}_{\alpha}\right\Vert
_{L^{\infty}(\Omega)}\right)  ^{2}}\left\Vert m^{-}\right\Vert _{L^{\infty
}(\Omega)}\leq\frac{1}{\lambda_{1}(m,I)},\label{i1}\\
\text{where\quad}\gamma_{b}:=\max\left\{  \left\Vert \overline{B}_{\alpha
}\right\Vert _{L^{1}(\alpha,x_{1})},\left\Vert \overline{B}_{\alpha
}\right\Vert _{L^{1}(x_{0},\beta)}\right\}  ,\nonumber
\end{gather}
then there \textit{exists }$u\in W^{2,2}(\Omega)$\textit{\ solution of
}(\ref{prob})\textit{.}\newline(ii) If $c\equiv0$ and it holds that%
\begin{gather}
\left(  1-p\right)  \mathcal{M}<\frac{1}{\lambda_{1}(m,I)}\text{,\qquad
where}\label{i2}\\
\mathcal{M}:=\max\left\{  \int_{x_{0}}^{\beta}\overline{B}_{\alpha}\left(
x\right)  \left\Vert m^{-}\underline{B}_{\alpha}\right\Vert _{L^{1}\left(
x,\beta\right)  }dx,\int_{\alpha}^{x_{1}}\overline{B}_{\alpha}\left(
x\right)  \left\Vert m^{-}\underline{B}_{\alpha}\right\Vert _{L^{1}\left(
\alpha,x\right)  }dx\right\}  ,\nonumber
\end{gather}
then there \textit{exists }$u\in W^{2,2}(\Omega)$\textit{\ solution of
}(\ref{prob})\textit{.}\newline
\end{theorem}

\textit{Proof}. The proof follows the lines of the proof of Theorem \ref{aa}
and hence we omit the details. Let us prove (i). We take $u_{2}$ as in the
aforementioned theorem, and we choose $\tau$ such that
\[
\frac{\left(  \gamma_{b}\left\Vert \underline{B}_{\alpha}\right\Vert
_{L^{\infty}(\Omega)}\right)  ^{2}}{C_{p}-\left\Vert c\right\Vert _{L^{\infty
}(\Omega)}\left(  \gamma_{b}\left\Vert \underline{B}_{\alpha}\right\Vert
_{L^{\infty}(\Omega)}\right)  ^{2}}\left\Vert m^{-}\right\Vert _{L^{\infty
}(\Omega)}\leq\frac{1}{\tau}\leq\frac{1}{\lambda_{1}(m,I)}.
\]
Let $x\in\left[  \alpha,x_{1}\right]  $ and define
\begin{gather*}
u_{1}\left(  x\right)  :=\left(  \sigma\int_{\alpha}^{x}\overline{B}_{\alpha
}\left(  y\right)  dy\right)  ^{k},\qquad\text{where}\\
\sigma:=\left[  \frac{\left\Vert \underline{B}_{\alpha}\right\Vert
_{L^{\infty}(\Omega)}^{2}\left(  \tau\left\Vert m^{-}\right\Vert _{L^{\infty
}(\Omega)}+\left\Vert c\right\Vert _{L^{\infty}(\Omega)}\right)  }{C_{p}%
}\right]  ^{1/2},\qquad k:=\frac{2}{1-p}.
\end{gather*}
We have that $u_{1}(\alpha)=0$, $u_{1}>0$ in $(\alpha,x_{1})$ and that $u_{1}$
is increasing. Moreover, after some computations one can check that
$\left\Vert u_{1}\right\Vert _{\infty}\leq1$ and
\begin{gather*}
-\left(  \underline{B}_{\alpha}\left(  x\right)  u_{1}^{\prime}\left(
x\right)  \right)  ^{\prime}=k\left(  k-1\right)  \sigma^{2}\left(  \sigma
\int_{\alpha}^{x}\overline{B}_{\alpha}\left(  y\right)  dy\right)
^{k-2}\overline{B}_{\alpha}\left(  x\right)  \leq\\
-\left\Vert \underline{B}_{\alpha}\right\Vert _{L^{\infty}(\Omega)}\left(
\tau\left\Vert m^{-}\right\Vert _{L^{\infty}(\Omega)}+\left\Vert c\right\Vert
_{L^{\infty}(\Omega)}\right)  \left(  \sigma\int_{\alpha}^{x}\overline
{B}_{\alpha}\left(  y\right)  dy\right)  ^{kp}\leq\\
\underline{B}_{\alpha}\left(  \tau m-c\right)  u_{1}^{p}\leq\underline
{B}_{\alpha}\left(  \tau mu_{1}^{p}-cu_{1}\right)  ,
\end{gather*}
that is, $Lu_{1}\leq\tau mu_{1}^{p}$ in $(\alpha,x_{1})$. The existence of
$u_{3}$ follows similarly. Let us prove (ii). We pick $\tau$ satisfying%
\begin{equation}
\left(  1-p\right)  \mathcal{M}<\frac{1}{\tau}<\frac{1}{\lambda_{1}(m,I)}.
\label{uso}%
\end{equation}
For $x\in\left[  \alpha,x_{1}\right]  $ we define
\begin{gather*}
u_{1}\left(  x\right)  :=\left(  \sigma\int_{\alpha}^{x}\overline{B}_{\alpha
}\left(  y\right)  \left\Vert m^{-}\underline{B}_{\alpha}+\varepsilon
\right\Vert _{L^{1}\left(  \alpha,y\right)  }dy\right)  ^{k}\qquad
\text{where}\\
\sigma:=\tau\left(  1-p\right)  ,\qquad k:=\frac{1}{1-p},\qquad\varepsilon>0.
\end{gather*}
Taking $\varepsilon$ small enough and employing (\ref{uso}) one can see that
$\left\Vert u_{1}\right\Vert _{\infty}\leq1$. Also, a few computations yield
that%
\begin{gather*}
-\left(  \underline{B}_{\alpha}\left(  x\right)  u_{1}^{\prime}\left(
x\right)  \right)  ^{\prime}\leq\\
-k\sigma^{k}\left(  \int_{\alpha}^{x}\overline{B}_{\alpha}\left(  y\right)
\left\Vert m^{-}\underline{B}_{\alpha}+\varepsilon\right\Vert _{L^{1}\left(
\alpha,y\right)  }dy\right)  ^{k-1}\left(  m^{-}\left(  x\right)
\underline{B}_{\alpha}\left(  x\right)  +\varepsilon\right)  \leq\\
-\tau m^{-}\left(  x\right)  \underline{B}_{\alpha}\left(  x\right)  \left(
\sigma\int_{\alpha}^{x}\overline{B}_{\alpha}\left(  y\right)  \left\Vert
m^{-}\underline{B}_{\alpha}+\varepsilon\right\Vert _{L^{1}\left(
\alpha,y\right)  }dy\right)  ^{kp}\leq\tau\underline{B}_{\alpha}mu_{1}^{p}.
\end{gather*}
Since $u_{3}$ can be defined analogously, this concludes the proof of (ii).
$\blacksquare$

\begin{remark}
\label{ojito1}Let us note that the inequalities in (i) and (ii) are not
comparable because one involves the $L^{\infty}$-norm of $m^{-}$ and the
constant $C_{p}$, and the other one does not. $\blacksquare$
\end{remark}

\begin{remark}
\label{ojito2}(i) It can be verified that (\ref{seno}) is better than
(\ref{i1}) when $b\equiv0$ (noting that in this case $\underline{B}_{\alpha
}=\overline{B}_{\alpha}=1$ and $\gamma_{b}=\gamma$ ($\gamma$ as in the
statement of Theorem \ref{aa})). If also $c\equiv0$, (\ref{i1}) becomes
exactly (\ref{lap}), that is, the condition deduced from the aforementioned
theorem for the laplacian operator. \newline(ii) In the case $b\equiv0$,
(\ref{i2}) reads as
\begin{equation}
\left(  1-p\right)  \max\left\{  \int_{x_{0}}^{\beta}\left\Vert m^{-}%
\right\Vert _{L^{1}\left(  t,\beta\right)  }dt,\int_{\alpha}^{x_{1}}\left\Vert
m^{-}\right\Vert _{L^{1}\left(  \alpha,t\right)  }dt\right\}  <\frac
{1}{\lambda_{1}(m,I)} \label{b0}%
\end{equation}
which is substantially better than the condition stated in \cite{ultimo},
Theorem 2.1, for $L=-u^{\prime\prime}$. Also, (\ref{b0}) is clearly not
comparable (for the same reason as in the above remark) with the inequalities
that can deduced from Theorem \ref{aa} in the case $c\equiv0$ (i.e., as the
one included in Remark \ref{lapla}). $\blacksquare$
\end{remark}

\begin{corollary}
\label{puf}Let $K_{b}:=\int_{\alpha}^{\beta}\overline{B}_{\alpha}\left(
x\right)  \left\Vert \underline{B}_{\alpha}\right\Vert _{L^{2}(\alpha,x)}dx$.
If (\ref{i2}) holds with $m/\left(  K_{b}\left\Vert m^{+}\right\Vert
_{L^{2}(\alpha,\beta)}\right)  -c$ instead of $m$, then there \textit{exists
}$u\in W^{2,2}(\Omega)$\textit{\ solution of }(\ref{prob})\textit{.}
\end{corollary}

\textit{Proof}. Applying H\"{o}lder%
%TCIMACRO{\U{b4}}%
%BeginExpansion
\'{}%
%EndExpansion
s inequality in (\ref{iner}) we see that $\left\Vert u\right\Vert _{L^{\infty
}(\Omega)}^{1-p}\leq\left\Vert m^{+}\right\Vert _{L^{2}(\alpha,\beta)}K_{b}$
for any nonnegative subsolution of (\ref{prob}). Now, let $\tau:=1/\left(
K_{b}\left\Vert m^{+}\right\Vert _{L^{2}(\alpha,\beta)}\right)  $, and let $u$
be the solution of (\ref{prob}) with $\tau m-c$ in place of $m$ provided by
Theorem \ref{bien} (ii). It follows that $\left\Vert u\right\Vert _{\infty
}\leq1$ and thus
\[
-u^{\prime\prime}+bu^{\prime}=\left(  \tau m-c\right)  u^{p}\leq\tau
mu^{p}-cu
\]
and recalling once again Remark \ref{homsup} the corollary follows.
$\blacksquare$

\begin{remark}
\label{lm}(i) Given \textit{any operator }$L$\textit{\ and any }$m\in
L^{2}\left(  \Omega\right)  $\textit{\ with }$0\not \equiv m\geq0$\textit{\ in
some }$I\subset\Omega$, let us note that the above corollary implies that
(\ref{prob}) has a solution if $p$ is sufficiently close to $1$.\newline(ii)
Given \textit{any operator }$L$\textit{\ and any }$m\in L^{2}\left(
\Omega\right)  $\textit{\ with }$m^{-}\in L^{\infty}\left(  \Omega\right)
$\textit{\ and }$0\not \equiv m\geq0$\textit{\ in some }$I\subset\Omega$, let
us observe that (\ref{i1}) says that (\ref{prob}) possesses a solution for
$\overline{m}:=m\chi_{\Omega-I}+km\chi_{I}$ if $k>0$ is large enough.
$\blacksquare$
\end{remark}

The next result provides the structure of the set of $p^{\prime}$s such that
(\ref{prob}) has a solution.

\begin{corollary}
\label{ppp}Let $m\in C\left(  M^{+}\right)  \cap L^{2}\left(  \Omega\right)  $
with $m^{+}\not \equiv 0$\textit{, }and let $\mathcal{P}$ be the set of
$p\in\left(  0,1\right)  $ such that (\ref{prob}) admits some solution $u\in
W^{2,2}\left(  \Omega\right)  $. Then $\mathcal{P}=\left(  0,1\right)  $ or
either $\mathcal{P}=\left(  p,1\right)  $ or $\mathcal{P}=\left[  p,1\right)
$ for some $p>0$.
\end{corollary}

\textit{Proof}. By Remark \ref{lm} (i) we have that $\mathcal{P}%
\not =\emptyset$. Let $p^{\ast}:=\inf\mathcal{P}$. If $\mathcal{P}%
\not =\left(  0,1\right)  $, Lemma \ref{qp} implies that $p^{\ast}>0$ and that
(\ref{prob}) has a solution for every $p>p^{\ast}$. Therefore, either
$\mathcal{P}=\left(  p^{\ast},1\right)  $ or $\mathcal{P}=\left[  p^{\ast
},1\right)  $. $\blacksquare$

\qquad

We write
\begin{gather}
I_{R}\left(  x_{0}\right)  :=\left(  x_{0}-R,x_{0}+R\right)  ,\label{jota}\\
\mathfrak{I}:=\left\{  I_{R}\left(  x_{0}\right)  \subset\Omega:m\leq0\text{
\textit{in }}I_{R}\left(  x_{0}\right)  \right\}  .\nonumber
\end{gather}

\begin{theorem}
\label{necee}Let $C_{p}$ and $\mathfrak{I}$ be given by (\ref{cp}) and
(\ref{jota}) respectively. Suppose there exists $u\in W^{2,2}\left(
\Omega\right)  $ solution of (\ref{prob}). Then
\begin{gather}
\sup_{I_{R}\left(  x_{0}\right)  \in\mathfrak{I}}\left[  \left[  \frac
{\gamma_{b,R}}{\left\Vert \overline{B}_{\alpha}\right\Vert _{L^{\infty}%
(I_{R}\left(  x_{0}\right)  )}}\right]  ^{2}\inf_{I_{R}\left(  x_{0}\right)
}m^{-}\right]  \leq C_{p}\int_{\alpha}^{\beta}\overline{B}_{\alpha}\left(
x\right)  \left\Vert m^{+}\underline{B}_{\alpha}\right\Vert _{L^{1}(\alpha
,x)}dx\text{,}\label{nec}\\
\text{\textit{where\quad}}\gamma_{b,R}:=\min\left\{  \int_{x_{0}}^{x_{0}%
+R}\overline{B}_{\alpha}\left(  y\right)  dy,\int_{x_{0}-R}^{x_{0}}%
\overline{B}_{\alpha}\left(  y\right)  dy\right\}  .\nonumber
\end{gather}
Let $M^{+}$ be given by (\ref{MM}). If also $c>0$ in $M^{+}$, then (\ref{nec})
must also hold with $C_{p}\sup_{x\in M^{+}}\frac{m^{+}\left(  x\right)
}{c\left(  x\right)  }$ in the right side of the inequality.
\end{theorem}

\textit{Proof}. We proceed by contradiction. Suppose (\ref{nec}) is not true
and let $I_{R}\left(  x_{0}\right)  \in\mathfrak{I}$ be such that
\begin{equation}
C_{p}\int_{\alpha}^{\beta}\overline{B}_{\alpha}\left(  x\right)  \left\Vert
m^{+}\underline{B}_{\alpha}\right\Vert _{L^{1}(\alpha,x)}dx\leq\left[
\frac{\gamma_{b,R}}{\left\Vert \overline{B}_{\alpha}\right\Vert _{L^{\infty
}(I_{R}\left(  x_{0}\right)  )}}\right]  ^{2}\inf_{I_{R}\left(  x_{0}\right)
}m^{-}\text{.} \label{berp}%
\end{equation}
For $x\in\overline{I}_{R}\left(  x_{0}\right)  $ we define a function $w\ $as
follows. If $x\in\left[  x_{0},x_{0}+R\right]  $ we set
\begin{gather*}
w\left(  x\right)  :=\left(  \sigma\int_{x_{0}}^{x}\overline{B}_{\alpha
}\left(  y\right)  dy\right)  ^{k},\qquad\text{where}\\
\sigma:=\left[  \frac{\inf_{I_{R}\left(  x_{0}\right)  }m^{-}}{C_{p}\left\Vert
\overline{B}_{\alpha}\right\Vert _{L^{\infty}(I_{R}\left(  x_{0}\right)
)}^{2}}\right]  ^{1/2},\qquad k:=\frac{2}{1-p},
\end{gather*}
and if $x\in\left[  x_{0}-R,x_{0}\right]  $ we set $w\left(  x\right)
:=\left(  \sigma\int_{x}^{x_{0}}\overline{B}_{\alpha}\left(  y\right)
dy\right)  ^{k}$ with $\sigma$ and $k$ as above. In $\left(  x_{0}%
,x_{0}+R\right)  $ we find that
\begin{gather*}
\left(  \underline{B}_{\alpha}w^{\prime}\right)  ^{\prime}-\underline
{B}_{\alpha}cw\leq k\left(  k-1\right)  \sigma^{2}\left(  \sigma\int_{\alpha
}^{x}\overline{B}_{\alpha}\left(  y\right)  dy\right)  ^{k-2}\overline
{B}_{\alpha}\leq\\
\frac{\inf_{I_{R}\left(  x_{0}\right)  }m^{-}}{\left\Vert \overline{B}%
_{\alpha}\right\Vert _{L^{\infty}(I_{R}\left(  x_{0}\right)  )}}\left(
\sigma\int_{\alpha}^{x}\overline{B}_{\alpha}\left(  y\right)  dy\right)
^{kp}\leq\underline{B}_{\alpha}m^{-}w^{p},
\end{gather*}
i.e., $Lw\geq-m^{-}w^{p}$, and the same is also valid in $\left(
x_{0}-R,x_{0}\right)  $.

Let $u$ be a solution of (\ref{prob}). We claim that $u\leq w$ in
$I_{R}\left(  x_{0}\right)  $. Indeed, if not, let $\mathcal{O}:=\left\{  x\in
I_{R}\left(  x_{0}\right)  :w\left(  x\right)  <u\left(  x\right)  \right\}
$. Since $Lu=-m^{-}u^{p}$ in $I_{R}\left(  x_{0}\right)  $, we have $L\left(
w-u\right)  \geq m^{-}\left(  u^{p}-w^{p}\right)  \geq0$ in $\mathcal{O}$. Let
$\overline{x}\in\partial\mathcal{O}$. Then $w\left(  \overline{x}\right)
=u\left(  \overline{x}\right)  $ or either $\overline{x}=x_{0}+R$ or
$\overline{x}=x_{0}-R$. If $\overline{x}=x_{0}+R$, by Lemma \ref{inerte} and
(\ref{berp}) we obtain%
\begin{gather*}
u\left(  \overline{x}\right)  ^{1-p}\leq\left\Vert u\right\Vert _{L^{\infty
}\left(  \Omega\right)  }^{1-p}\leq\int_{\alpha}^{\beta}\overline{B}_{\alpha
}\left(  x\right)  \left\Vert m^{+}\underline{B}_{\alpha}\right\Vert
_{L^{1}(\alpha,x)}dx\leq\\
\left[  \frac{\int_{x_{0}}^{x_{0}+R}\overline{B}_{\alpha}\left(  y\right)
dy}{\left\Vert \overline{B}_{\alpha}\right\Vert _{L^{\infty}(I_{R}\left(
x_{0}\right)  )}}\right]  ^{2}\frac{\inf_{I_{R}\left(  x_{0}\right)  }m^{-}%
}{C_{p}}=w\left(  \overline{x}\right)  ^{1-p}\text{,}%
\end{gather*}
and we arrive to the same inequality if\ $\overline{x}=x_{0}-R$. Therefore the
maximum principle says that $u\leq w$ in $\mathcal{O}$ which is not possible.
Thus, $u\leq w$ in $I_{R}\left(  x_{0}\right)  $; but $u>0$ in $\Omega$ and
$w\left(  x_{0}\right)  =0$. Contradiction.

To conclude the proof we note that the last statement of the theorem may be
derived as above applying Lemma \ref{dudu} instead of Lemma \ref{inerte}.
$\blacksquare$

\qquad

\begin{remark}
(i) It follows from the above theorem that given $b$, $m$, $p$ fixed, there
exists $0\leq c_{0}\in L^{\infty}(\Omega)$ such that for all $c\in L^{\infty
}(\Omega)$ with $c\geq c_{0}$ the problem (\ref{prob}) \textit{does not admit
a solution}. Note that given $L$, $m$, $p$ fixed with $0\not \equiv m\leq
0$\ in some $I\subset\Omega$, \textit{neither there is a solution} for
$\underline{m}:=m\chi_{\Omega-I}+km\chi_{I}$ if $k>0$ is large enough.\newline%
(ii) We observe that (\ref{nec}) always is true if $p$ is sufficiently close
to $1$. Let us mention that this must indeed occur by Remark \ref{lm}.
$\blacksquare$
\end{remark}

As a consequence of the previous theorems we derive an existence result for
problems of the form
\begin{equation}
\left\{
\begin{array}
[c]{ll}%
Lu=mf\left(  u\right)   & \text{en }\Omega\\
u>0 & \text{en }\Omega\\
u=0 & \text{en }\partial\Omega,
\end{array}
\right.  \label{f}%
\end{equation}
for certain continuous functions $f:\left[  0,\infty\right)  \rightarrow
\left[  0,\infty\right)  $. Now we state assumption

(H1) There exist $k_{1},k_{2}>0$ and $p\in\left(  0,1\right)  $ such that%
\begin{gather*}
k_{1}\xi^{p}\leq f\left(  \xi\right)  \leq k_{2}\xi^{p}\text{ for all }\xi
\in\left[  0,\underline{K}\right]  \text{, }\\
\text{where }\underline{K}:=\left[  k_{1}\int_{\alpha}^{\beta}\overline
{B}_{\alpha}\left(  x\right)  \left\Vert m^{+}\underline{B}_{\alpha
}\right\Vert _{L^{1}(\alpha,x)}dx\right]  ^{1/\left(  1-p\right)  },\\
\text{and }f\left(  \xi\right)  \leq k_{3}\xi^{q}\text{ for all }\xi\in\left[
\overline{K},\infty\right)  \text{ and some }\overline{K},k_{3}>0\text{ and
}q\in\left(  0,1\right)  .
\end{gather*}
Note that we make no monotonicity nor concavity assumptions on $f$.

\begin{corollary}
\label{fff}Let $f$ satisfy (H) and suppose (\ref{prob}) admits a solution with
$k_{1}m^{+}-k_{2}m^{-}$ instead of $m$. Then there exists $u\in W^{2,2}\left(
\Omega\right)  $ solution of (\ref{f}).
\end{corollary}

\textit{Proof}. Let $u\ $be the solution of (\ref{prob}) with $k_{1}%
m^{+}-k_{2}m^{-}$ in place of $m$. It follows from Lemma \ref{inerte} that
$\left\Vert u\right\Vert _{\infty}\leq\underline{K}$, and so from (H) we
deduce that
\[
Lu=\left(  k_{1}m^{+}-k_{2}m^{-}\right)  u^{p}\leq mf\left(  u\right)
\qquad\text{in }\Omega\text{.}%
\]
On the other side, let $\varphi>0$ be the solution of $L\varphi=m^{+}$ in
$\Omega$, $\varphi=0$ on $\partial\Omega$, and let $k\geq\max\left\{
\overline{K},(k_{3}\left(  \left\Vert \varphi\right\Vert _{\infty}+1\right)
^{q})^{1/(1-q)}\right\}  $. Recalling (H) and reasoning as in (\ref{ju}) we
see that
\[
L(k(\varphi+1))\geq km^{+}\geq k_{3}(k(\varphi+1))^{q}m^{+}\geq mf(k(\varphi
+1))\text{\qquad in }\Omega
\]
and the corollary is proved. $\blacksquare$

\end{document}